\documentclass[preprint,nonatbib]{elsarticle}

\usepackage{lineno,hyperref}
\usepackage{amsmath}
\usepackage{amsthm,amssymb}
\usepackage{verbatim}
\usepackage{mathptmx}      % use Times fonts if available on your TeX system
\usepackage{latexsym}
\usepackage{marvosym}

\journal{Journal of \LaTeX\ Templates}

\newtheorem{thm}{Theorem}

\newdefinition{rmk}{Remark}
\newproof{pf}{Proof}
\newproof{pot}{Proof of Theorem \ref{thm2}}
\newcommand{\SP}{\text{supp}}
\newcommand{\MYdef}{\mathrel{\stackrel{\rm def}=}}

\begin{document}
\begin{frontmatter}

\title{Constructive description of H{\"o}lder-like classes on an arc in $\mathbb{R}^3$ by means of harmonic functions}%\tnoteref{mytitlenote}}

\author[rvt]{Tatyana~A.~Alexeeva\corref{cor1}}
\ead{tatyanalexeeva@gmail.com}

\author[rvt,focal]{Nikolay~A.~Shirokov}
\ead{nikolai.shirokov@gmail.com}

\cortext[cor1]{Corresponding author}
\address[rvt]{Department of Applied Mathematics and Business Informatics, 

National Research University Higher School of Economics, 

3A Kantemirovskaya Ul., St.~Petersburg, 194100, Russia}

%\fntext[myfootnote]{Department of Applied Mathematics and Business Informatics, National Research University Higher School of Economics, 3A Kantemirovskaya Ul., St.~Petersburg, 194100, Russia.}
%\address{Department of Applied Mathematics and Business Informatics, National Research University Higher School of Economics, 3A Kantemirovskaya Ul., St.~Petersburg, 194100, Russia}

\address[focal]{Department of Mathematical Analysis, Faculty of Mathematics and Mechanics, 

St.~Petersburg State University, 28 Universitetsky prospekt, Peterhof, St.~Petersburg, 198504, Russia}

%\fntext[myfootnote]{National Research University Higher School of Economics, 3A Kantemirovskaya Ul., St.~Petersburg, 194100, Russia.}
%\fntext[myfootnote]{$^2$ Department of Mathematical Analysis, Faculty of Mathematics and Mechanics, St.~Petersburg State University, 28 Universitetsky prospekt, Peterhof, St.~Petersburg, 198504, Russia}

%\author[mymainaddress,mysecondaryaddress]{$^1$ Department of Applied Mathematics and Business Informatics, National Research University Higher School of Economics, 3A Kantemirovskaya Ul., St.~Petersburg, 194100, Russia}


\fntext[myfootnote]{The second author was supported by the RFBR, grant 17-01-00607}

\begin{abstract}
We give a constructive description of H{\"o}lder-like classes of functions on  chord-arc curves in $\mathbb{R}^3$ in terms  of a rate of approximation by harmonic functions in shrinking neighborhoods of those curve.
\end{abstract}

\begin{keyword}
Constructive description\sep H{\"o}lder classes \sep Harmonic functions \sep Chord-arc curves
\MSC[2010] 41A30\sep  41A27
\end{keyword}

\end{frontmatter}

%\linenumbers

\section{Introduction}

The   constructive description of classes of functions in terms of a rate of approximation   by functions taken from  specific sets  (polynomials, rational functions, entire functions, etc.) was initiated by  D.~Jackson and S.~N.~Bernstein at the beginning of the 20th century. Nevertheless, a natural   problem of a constructive description of H{\"o}lder classes on a segment in terms of a rate of approximation by algebraic polynomials was solved only in 1956 (see \cite{Dzyadyk56}, \cite[ch.~7]{Dzyadyk08}).
Since then problems concerning  constructive description of classes of functions defined on domains in the complex plane have played a central role in   approximation theory. Many authors were involved in the following problem: let $G$ be a Jordan region in the complex plane $\mathbb{C}$, and let $H(G)$ be a class of functions $f$ analytic in the interior $\stackrel{\circ}{G}$ of $G$ and continuous (or smooth)  on the closure of $G$. What is the scale of approximation of functions from $H(G)$ by algebraic polynomials which makes it possible to find  the rate of smoothness of relevant functions? V.~K.~Dzyadyk (\cite{Dzyadyk59,Dzyadyk62,Dzyadyk63,Dzyadyk_Nikol63}) introduced a special type of weights $\rho_{1/n}(z)$ on the boundary $\Gamma$ of   $G$ such that the condition that $f$ is   analytic  in  $\stackrel{\circ}{G}$ and satisfies the $\alpha$-H{\"o}lder condition, $\alpha > 0$, $\alpha\notin \mathbb{N}$, is equivalent to the possibility of approximating $f$ by polynomials $P(z)$ of degree $\le n$ with the property 

\begin{equation*}
\ |f(z)-P_n(z)|\le C_f \rho_{1/n}^{\alpha}(z),\, z\in\Gamma . \eqno(\star) 
\end{equation*}

So,  for various regions in $\mathbb{C}$, the weights $\rho_{1/n}^{\alpha}(z)$ were a successful scale for a constructive description of the above-mentioned classes of functions. The main problem in that direction was to weaken the assumptions concerning the boundary $\Gamma$. The results progressed from a piecewise smooth in some sense \cite{Dzyadyk66,Lebedev71,Shirokov72} to a chord-arc \cite{Dzyadyk77} and finally to a quasiconformal property of a Jordan curve $\Gamma$ \cite{Belyi77}.

If turned out that if a function $f$ can be approximated by polynomials $P_n(z)$ of degree $\le n$ as in ($\star$), then $f$ is analytic in $\stackrel{\circ}{G}$ and satisfies the $\alpha$-H{\"o}lder condition for any Jordan domain $G$ \cite{Lebedev66,LebedevT70}.

In the case where the boundary $\Gamma = \partial G$ has  cusps,  the polynomial approximation with the rate $const\cdot \rho_{1/n}^{\alpha}(z)$ is appropriate not for all functions satisfying the $\alpha$-H{\"o}lder condition \cite{Shirokov01,Andrievskii85}. Consequently,  in the case of an arbitrary Jordan region, the scale $\rho_{1/n}^{\alpha}(z)$ is not suited for  constructive description of the $\alpha$-H{\"o}lder classes by means of complex polynomials.  This circumstance stimulated the introduction of  a modified scale $\rho_{1/n}^{\star\alpha}(z)$ \cite{Shirokov01,Andrievskii85,Andrievskii94}, which was used for constructive description of $\alpha$-H{\"o}lder classes in Jordan domains with non-empty interior.

In the case where the interior of $G$ is empty, i.e.,  if $ G=\Gamma$, the problem of a constructive description of H{\"o}lder (or H{\"o}lder-like) classes of functions defined on $\Gamma$ by means of their approximation by polynomials turned out to be more intricate. For example, if  $G=\Gamma_{\beta}\stackrel{def}{=} [-1,\, 0]\cup [0,\, e^{\beta}]$, $0 <\beta <\pi$, then a simple combination of $\rho_{1/n}(z)$ and $\rho_{1/n}^{\star}(z)$ cannot provide a constructive description of the $\alpha$-H{\"o}lder class \cite{Andrievskii94}. Even in the case of $\Gamma_{\beta}$, the answer is obtained with the help of a Cantor-like construction of a scale using both scales $\rho_{1/n}^{\alpha}(z)$ and $\rho_{1/n}^{\star\alpha}(z)$ \cite{Shirokov77}.

V.~V.~Andrievskii \cite{Andrievskii94} found an alternative approach to the problem of a constructive approximation of functional classes on   Jordan arcs. He used a uniform approximation of a function $f$ defined on a Jordan arc $L$ by polynomials $P_n$ along with uniform estimates of $P'_n(z)$ in a neighbourhood of $L$. We notice that  harmonic polynomials can also be used for a constructive description of H{\"o}lder-like classes of functions on continuums in $\mathbb{C}$ (V.~V.~Andrievskii, \cite{Andrievskii86,Andrievskii88}).

 We emphasize that all above-mentioned constructions of the scales $\rho_{1/n}^{\alpha}(z)$ and $\rho_{1/n}^{\star\alpha}(z)$ and constructive descriptions of H{\"o}lder classes on curves are applicable only for plane curve since each of these constructions uses  a conformal mapping of the complement $\mathbb{C}\setminus \overline{G}$ onto the exterior of the unit disc $\mathbb{D}$. However, the same problems can be considered for H{\"o}lder spaces on curves lying in arbitrary $\mathbb R^n$ or $\mathbb C^n$. 
 
 In the present paper, we obtain a constructive description of $H^{\alpha} (L)$ for   chord-arc curves $L$ lying in $\mathbb{R}^3$. 
As approximating functions, we use harmonic functions with certain estimates of their gradients in neighborhoods of a curve. The neighborhoods are connected with the rate of approximation -- they shrink when the approximation is getting better.  

The paper is organized as follows. In Section~2 we introduce notation and state our main results. Section~3 contains the proof of Theorem~4. Section~4 contains the proof of Theorem~5. Section~5 is concerned with the proof of the main result of the paper -- Theorem~1. Section~6 is devoted to the proof of Theorem~2.

\section{Main results}

%page1
We say that a non-closed Jordan curve $L\subset \mathbb{R}^3$ has a chord-arc property (or is a chord-arc curve) if there exists a constant $C=C(L)$ such that the length of the subarc $L$ between $M_1$ and $M_2$ does not exceed $C\cdot \Vert M_1 M_2 \Vert$ for all points $M_1$,  $M_2\in L$.
We denote by $B_{r}(M)$ an open ball in $\mathbb{R}^3$ with center $M$ and   radius $r$ and put $\Omega_{\delta}(L)=\bigcup_{M\in L}  B_{\delta}(M)$. Let $H^{\omega}(L)$ be the space of all complex-valued functions $f$ that are defined on $L$ and satisfy  the condition $|f(M_2)-f(M_1)|\leq C_{f}\omega \big ( \Vert M_1 M_2\Vert\big )$, where $\omega$ is a modulus of continuity with the property

%\cite{RefB} and \cite{RefJ}.

%(1)
\begin{equation}
\int_{0}^{x}\frac{\omega (t)}{t} dt \le C' \omega (x),\; x\int_{t}^{\infty}\frac{\omega (t)}{t^2} dt\le C'' \omega (x)
\label{f1}
\end{equation}
\noindent
(here and below we denote by $C$, $C'$, $C_1$, $\ldots$ various constants). One of our two  main results in the present paper is the following
theorem.

%Theorem_1
\begin{thm}
Assume that $L$ is a bounded non-closed chord-arc curve and $f\in H^{\omega}\big ( L \big )$. Then there exist constants $C_1=C_1 (f, L)$ and $C_2=C_2 (f, L)$ such that for every $\delta >0$ there exists a function $\upsilon_\delta$ harmonic in $\Omega_{\delta} (L)$ such that

%(2)
\begin{equation}
\ |\upsilon_\delta(M)-f(M)|\leq C_1\omega(\delta), \; M\in L
\label{f2}
\end{equation}
%(3)
\begin{equation}
\ |\nabla \upsilon_\delta(M)|\leq C_2 \frac{\omega(\delta)}{\delta}, \; M\in \Omega_{\delta}(L)\backslash\Omega_{\frac{\delta}{2}}(L)
\label{f3}
\end{equation}
\end{thm}
%page2'
Theorem~1 may be called "a direct theorem" of approximation like many similar statements concerning approximation by polynomials, rational functions, etc. The "converse theorem" to   Theorem~1 is also valid: if we take a  unit vector $\vec{\ell}$, then (3) implies that

%(3')
$$
\ |\upsilon_{\delta\ell}' (M)|\leq C_2\frac{\omega(\delta)}{\delta}, \, M\in \Omega_{\delta}(L)\backslash\Omega_{\frac{\delta}{2}}(L),\eqno(3')
$$

\noindent
and the maximum principle for a function $\upsilon_{\delta\ell}'$ harmonic in $\Omega_{\delta}(L)$ guarantees that   estimate (${3}^{\prime}$) is valid for $M\in\Omega_{\delta}$; this gives the estimate

%(3'')
$$
\ |\nabla \upsilon_{\delta}(M)|\leq C'_2 \frac{\omega(\delta)}{\delta}, \; M\in \Omega_{\delta}(L). \eqno(3'')
$$
Further, if $M_1, \, M_2\in L$ and $\|M_1M_2\| \leq \frac{\delta}{2}$, then the segment $\overrightarrow{M_{1}M_{2}}$ lies in $\Omega_{\delta}(L)$. Putting $\overline{\nu}=\frac{1}{\|M_1M_2\|}\cdot \; \overrightarrow{M_{1}M_{2}}$, we get

$$
\ f(M_2)-f(M_1)=\left( f(M_2)-\upsilon_{\delta}(M_2)\right)- \left(f(M_1)-\upsilon_{\delta}(M_1)\right)-
$$
%(3^0)
$$
\ -\int\limits_0^1 \upsilon_{\delta\overline{\nu}}' (M_1 + \overline{\nu} \|M_1M_2\|t) dt. \eqno(3^{\circ})
$$

\noindent
So, if we suppose that a function $f$ can be approximated by functions $\upsilon_{\delta}$ as in $\eqref{f2}$ and $\eqref{f3}$, then ($3^{\circ}$) and ($3'$) imply that $f\in \mathbf{H}^{\omega}(L)$.
The constructive description of the space $\mathbf{H}^{\omega}(L)$ in terms of functions $\upsilon_{\delta}$ harmonic in  $\Omega_{\delta}(L)$ is in a sense strict. This is the assertion of the second main result.

%Theorem_2
\begin{thm}
Let $1>\delta_{k}>0$, $\delta_{k}\to 0$, $\delta_{k}$ be monotonically decreasing, $\ell_{k}\to +\infty$, and let the modulus of continuity $\omega (t)$ satisfy conditions $\eqref{f1}$. Then there exists a function $f_0 \in \mathbf{H}^{\omega}\left([A_0,\, B_0]\right)$, where $A_0=(-1,\, 0,\, 0)$ and $B_0=(1,\, 0,\, 0)$, that cannot be approximated  by functions $V_{k}$ harmonic in the domains $\Omega_{\ell_{k}\delta_{k}}\left([A_0,\, B_0]\right)$ in the following way:

%(4)
\begin{equation}
\left|V_{k}(M)-f_0(M)\right|\leq C'_{1}\omega (\delta_{k}), \, M\in [A_0,\, B_0]
\end{equation}
\noindent
if the functions $V_{k}$ satisfy the condition

%(5)
\begin{equation}
\left|\nabla V_{k}(M)\right|\leq C'_{2}\frac{\omega (\delta_{k})}{\delta_{k}}, \, M\in\Omega_{\ell_{k}\delta_{k}}\left([A_0,\, B_0]\right)
\end{equation}
\end{thm}

%page3
The proof of   Theorem~1 depends on a special type of an extension of a function $f$ from the curve $L$ to the entire space $\mathbb{R}^{3}$; we call this extension \emph{pseudoharmonic} by analogy with the widely-used pseudoanalytic extension due to E.~M.~Dyn'kin \cite{Dynkin77,Dynkin93}.

%Theorem_3a=4
\begin{thm}
{\it Let  $f\in H^{\omega}\big ( L \big )$, where $\omega$ is a modulus of continuity satis\-fying assumption (1). Let  $O$ be the origin of $\mathbb{R}^{3}$. Then there is a function $f_0 \in C\left(\mathbb{R}^{3}\right)$ such that $f_0 |_{L}=f$, $f_0\in C^{2}\left(\mathbb{R}^{3}\setminus L\right)$, and}
%(6)
\begin{equation}
\left|\nabla f_{0}(M)\right|=o(dist^{-1}(M,L)), \; o\; \text{is uniform on}\; \mathbb{R}^{3},
\end{equation}
%(7)
\begin{equation}
f_{0}(M)\equiv 0, \;\text{for}\; \|\overrightarrow{OM}\| \geq R_{0}, \; \text{and}\; L\subset B_{R_{0}}(O)
\end{equation}
%(8)
\begin{equation}
\left|\Delta f_{0}(M)\right|\leq C_0\frac{\omega\left( dist(M,L)\right)}{dist^2\left(M,L\right)}
\end{equation}
\end{thm}
In what follows, we call an extension $f_0$ of a function $f$ a pseudoharmonic extension of $f$.

%Theorem_3b=5
\begin{thm}
{\it Assume that  a function $f\in C(L)$ has a pseudoanalytic extension satisfying conditions (6), (7), and (8). Then $f\in H^{\omega}\left( L \right)$}.
\end{thm}

Theorems~4 and 5 are exactly analogous to the theorems of E.~M.~Dyn'kin concerning pseudoanalytic extensions of   functions defined on domains in $\mathbb{C}$ \cite{Dynkin77,Dynkin93}.

\subsection{Proof of Theorem~4}
\label{sec:2}
%{\bf 3. Proof of Theorem~3a}

%page4
We begin with the proof of Theorem~4. Let $A$ be one of endpoints of the curve $L$ and let $B$ be the another one. In the sequel, we denote by $\ell (M_1,\, M_2)$   the length of the arc of $L$ with the endpoints $M_1$ and $M_2$. Let $\ell (A,\, B)=\Lambda$. We subdivide $L$ into $2^{n}$ arcs of equal length by the points $M_{kn}$, $0\leq k\leq 2^{n}$, $M_{0n}=A$, $M_{2^n,n}=B$, where the index  $k$ increases as the points $M_{kn}$ move in the direction from $A$ to $B$. The chord-arc property of $L$ implies the inequality

$$
\|\overrightarrow{M_{kn},\,M_{k+1,n}}\| \geq \frac{1}{C_0} \ell\left(M_{kn},\, M_{k+1,n}\right)=\frac{1}{C_0}\cdot 2^{-n} \Lambda \stackrel{\rm def}=\frac{1}{C_0} \Lambda_{n}.
$$

\noindent
We put

%(9)
\begin{equation}
\widetilde{\Omega}_{n}\stackrel{\rm def}=\bigcup_{k=0}^{2^n} \overline{B}_{2\Lambda_{n}}(M_{kn}),
\end{equation}

%(10)
\begin{equation}
\Omega_{n}\stackrel{\rm def}=\widetilde{\Omega}_{n}\setminus \overline{\widetilde{\Omega}}_{n+1}.
\end{equation}

\noindent
For $M\in \Omega_{n}$ we have the estimates

%(11)
\begin{equation}
\frac{1}{2}\Lambda_{n}\leq dist\left(M,\,L\right)\leq 2\Lambda_{n}.
\end{equation}

\noindent
Let

%(12)
\begin{equation}
\omega_{0n}=B_{2\Lambda{n}}\left(M_{0n}\right)\cap{\Omega_{n}},
\end{equation}

%(13)
\begin{equation}
\omega_{kn}=\left(B_{2\Lambda_{n}}\left(M_{kn}\right)\cap{\Omega_{n}}\right)\setminus \bigcup_{\nu=0}^{k-1} B_{2\Lambda_{n}}\left(M_{\nu n}\right),\; 1\leq k\leq 2^{n}
\end{equation}
($\omega_{kn}$ may be empty for some $k$ and $n$). We define the function $g$ as follows:

%(14)
\begin{equation}
g(M)=\left\{
\begin{array}{rl}
f\left(M_{kn}\right),\, M\in\omega_{kn} \\
0, \, M\in \mathbb{R}^{3}\setminus\bigcup_{n=0}^{\infty} \widetilde{\Omega}_{n}
\end{array}\right.
\end{equation}

%page5
Let $d(M)=dist\left(M,\,L\right)$, $M\in \mathbb{R}^{3}\setminus L$ and $B_{\star}\left(M\right)=\overline{B}_{\frac{1}{8}d(M)}\left(M\right).$
\noindent
We need to control the distance $\|\overrightarrow{M_{kn},\, M_{k_1,\, n_1}}\|$ in the case where $M\in \omega_{kn}$, $M_1\in B_{\star}(M)\cap\omega_{k_1,\, n_1}$. We have $$2\Lambda_{n_1}\ge d(M_1)\ge d(M)-\|\overrightarrow{M M_1}\|\ge\frac{1}{2}\Lambda_{n}-\frac{1}{8} d(M)\ge \frac{1}{2}\Lambda_{n}-2\cdot\frac{1}{8}\Lambda_{n}=\frac{1}{4}\Lambda_{n},$$ from which we obtain $8\Lambda_{n_{1}}\ge\Lambda_{n}$, $-n_{1}+3\ge -n$, and $n_{1}\leq n+3$.

\noindent
Then we observe that
$$
\frac{1}{2}\Lambda_{n_1}\leq d(M_1)\leq d(M) +\|\overrightarrow{M M_1}\|\leq
$$
$$
\leq 2\Lambda_{n}+\frac{1}{8} d(M)\leq 2\Lambda_{n} + 2\cdot\frac{1}{8}\Lambda_{n} = 2\frac{1}{4}\Lambda_{n} <4\Lambda_{n} ,
$$

\noindent
hence $\Lambda_{n_1} < 8\Lambda_{n}$, $-n_1 < -n+3$, and $n_1\ge n-2$.

\noindent
Let $N, N_1\in L$  be such  that $\|\overrightarrow{MN}\| = d(M)$, $\|\overrightarrow{M_1 N_1}\| = d(M_1)$. Since $\|\overrightarrow{N M_{kn}}\| \leq 4\Lambda_{n}$, $\|\overrightarrow{N_1 M_{k_1 n_1}}\| \leq 4\Lambda_{n_1}$,
and
$$
\|N N_1\|\leq \|NM\|+\|M M_1\| + \|M_1 N_1\|\leq 2\Lambda_{n}+\frac{1}{8}\Lambda_{n} +2\Lambda_{n_1} \leq
$$
$$
\leq 2\Lambda_{n}+\frac{1}{8}\Lambda_{n} +2\cdot 8\Lambda_{n}=18\frac{1}{8}\Lambda_{n} ,
$$

\noindent
we have the estimates
%(15)
\begin{equation}
\begin{split}
&\|M_{kn} M_{k_1 n_1}\|\leq \|M_{kn} N\| +\|N N_1\| +\|N_1 M_{k_1 n_1}\|\leq  \\
& \leq 4\Lambda_{n} + 18\frac{1}{8}\Lambda_{n} +4\Lambda_{n_1}\leq \left(4+18\frac{1}{8} +4\cdot 8\right) < 55\Lambda_{n}
\end{split}
\end{equation}

\noindent
Inequality (15) and assumption ($8'$) imply the inequalities

%(16)
\begin{equation}
\left| f\left(M_{kn}\right) - f\left(M_{k_1 n_1}\right)\right|\leq \omega\left(55\Lambda_{n}\right)\leq C\omega\left(\Lambda_{n}\right) .
\end{equation}
As a consequence of (16) and (14) we get the inequality
%(17)
\begin{equation}
\left| g\left(M_1\right) - g\left(M\right)\right|\leq C\omega\left(d(M)\right).
\end{equation}
valid for all $M_{1}\in B_{\star}(M)$.
We define

%(18)
\begin{equation}
g_{1}(M)=\frac{1}{\left| B_{\star}(M)\right|}\int\limits_{B_{\star}(M)} g\left( M_1\right) d{m_3} ,
\end{equation}
where $\left| B_{\star}(M)\right|$  is the volume of the ball $B_{\star}(M)$ and $m_3$ is the 3-dimensional Lebesgue measure. Due to (18) and (17) we see that $g_1\in C\left(\mathbb{R}^{3}\setminus L\right)$ and

%(19)
\begin{eqnarray}
&\left| g_1\left(M\right) - g\left(M\right)\right| =\left| g_1\left(M\right) - f\left(M_{kn}\right)\right|= \nonumber \\
&=\left|\frac{1}{\left| B_{\star}(M)\right|}\int\limits_{ B_{\star}(M)} g(M_1)\, d{m_3}\left( M_1\right) - \frac{1}{\left| B_{\star}(M)\right|}\int\limits_{ B_{\star}(M)} g(M)\, d{m_3}\left( M_1\right)\right|\leq \\
&\leq C\omega\left(d(M)\right) \nonumber .
\end{eqnarray}
The definition (14) and  estimate (19) imply that $g_{1}(M)\rightarrow f(M_{\star})$ as $M\rightarrow M_{\star}$, $M_{\star}\in L$. Hence  the function $g_1$ is continuous on $\mathbb{R}^{3}$ and vanishes outside a certain ball.

Now we construct   a characteristic $d_{0}(M)$ that is commeasurable with $d(M)$ but is $C^2\left(\mathbb{R}^{3}\setminus L\right)$-smooth in contrast to $d(M)$, which is usually only Lip1 on $\mathbb{R}^{3}\setminus L$. Let $\sum_{n}^{} = \{ M\in \mathbb{R}^{3}\setminus L :\, 2^{n-1} < d(M)\leq 2^{n}\}$, $n\in \mathbb{Z}$. Since
$$
\left| d(M_2) - d(M_1)\right|\leq \|M_1 M_2\|,\; M_1, M_2 \in \mathbb{R}^{3}\setminus L ,
$$
\noindent
the balls $B_{r_1}\left( M_1\right)$ and $B_{r_2}\left( M_2\right)$ are disjoint if $r_1 < \frac{1}{4} d(M_1)$, $r_2 < \frac{1}{4} d(M_1)$, and $d(M_2)\ge 2d(M_1)$. Due to this observation,  the following functions are well defined:

%page7
%(20)
\begin{equation}
d_{1}(M)=2^{n-1}, \; M\in\sum_{n}, \; n\in\mathbb{Z}
\end{equation}

%(21)
\begin{equation}
d_{2}(M)=\frac{1}{\left| B_{\frac{1}{8}\cdot 2^{n-1}}(M)\right|}\int\limits_{ B_{\frac{1}{8}\cdot 2^{n-1}}(M)} d_{1}(\widetilde{M})\, dm_{3}(\widetilde{M}),
\end{equation}

\noindent
if $2^{n-1}\cdot\sqrt{2}<d(M)\leq \sqrt{2}\cdot 2^{n}=\frac{1}{\sqrt{2}}\cdot 2^{n+1}$.
We observe that definitions (20) and (21) imply the estimate $\|\operatorname{grad}d_{2}(M)\|\leq C$. Finally, we put

%(22)
\begin{equation}
d_{0}(M)=\frac{1}{\left| B_{\frac{1}{8}\cdot 2^{n-1}}(M)\right|}\int\limits_{ B_{\frac{1}{8}\cdot 2^{n-1}}(M)} d_{2}(\widetilde{M})\, dm_{3}(\widetilde{M}),
\end{equation}

\noindent
if $2^{n-1}\cdot\sqrt{2}<d(M)\leq \sqrt{2}\cdot 2^{n}$.

Equation (22) gives the required function $d_0$. We have the following estimates:

%(23)
\begin{equation}
d_{0}(M)\asymp d(M), \; \|\operatorname{grad}d_{0}(M)\|\leq C
\end{equation}

\noindent
and

%(24)
\begin{equation}
\|\operatorname{grad}^{2}d_{0}(M)\|\leq \frac{C}{d(M)},
\end{equation}

\noindent
which follow from  (22). Indeed, if $\bar{\lambda}$, $\bar{\mu}$ are arbitrary unit vectors, then (22) implies

$$
{d'_{0}}_{\bar{\lambda}}(M)=\frac{1}{\left| B_{\frac{1}{8}\cdot 2^{n-1}}(M)\right|}\int\limits_{B_{\frac{1}{8}\cdot 2^{n-1}}(M)} {d'_{2}}_{\bar{\lambda}}(\widetilde{M})\, dm_{3}(\widetilde{M}),
$$

\noindent
which gives (23), and if $\bar{\nu}(\widetilde{M})$ is the outer unit normal to the sphere $S_{\frac{1}{8}\cdot 2^{n-1}}(M)$ at the point $\widetilde{M}$,
then
%(25)
\begin{equation}
{d''_{0}}_{\bar{\lambda}\bar{\mu}}(M)=\frac{1}{\left| B_{\frac{1}{8}\cdot 2^{n-1}}(M)\right|}\int\limits_{B_{\frac{1}{8}\cdot 2^{n-1}}(M)}\left(\bar{\mu}, \,\bar{\nu}(M) \right) {d'_{2}}_{\bar{\lambda}}(\widetilde{M})\, dS(\widetilde{M}),
\end{equation}
where
$dS(\widetilde{M})$ denotes the Lebesgue measure on $S_{\frac{1}{8}\cdot 2^{n-1}}(\widetilde{M})$;   estimate (24) follows from (25). Let us notice that $d_{1}(M)\leq d(M)$, and, for $\widetilde{M}\in B_{\frac{1}{8}\cdot 2^{n-1}}(M)$, we also have $d_{1}(\widetilde{M})\preccurlyeq d(M)$, hence $d_{2}(M)\preccurlyeq d(M)$. Moreover, (22) implies that $d_{0}(M)\preccurlyeq d(M)$. Finally, we define

%(26)
\begin{equation}
g_{2}(M)=\frac{1}{\left| B_{\frac{1}{8} d_{0}(M)}(M)\right|}\int\limits_{B_{\frac{1}{8} d_{0}(M)}(M)} g_{1}(\widetilde{M})\, dm_{3}(\widetilde{M}),
\end{equation}

%(27)
\begin{equation}
g_{0}(M)=\frac{1}{\left| B_{\frac{1}{8} d_{0}(M)}(M)\right|}\int\limits_{B_{\frac{1}{8} d_{0}(M)}(M)} g_{2}(\widetilde{M})\, dm_{3}(\widetilde{M}),
\end{equation}

\noindent
We notice  that definitions (20)--(22) imply the inequalities $$d_{1}(M)\ge\frac{1}{2} d(M),  \ d_{2}(M)\ge\frac{1}{2} d(M), \ d_{0}(M)\ge\frac{1}{2} d(M).$$
Let $B^{\star}(M)=B_{\frac{1}{8} d_{0}(M)}(M)$ and $r^{\star}(M)=\frac{1}{8} d_{0}(M)$. Using these estimates in the same  way as in (19), we get the estimates

%(28)
\begin{equation}
\left| g_2(M)-g(M)\right|\leq C\omega (d(M))
\end{equation}

\noindent
and

%(29)
\begin{equation}
\left| g_0(M)-g(M)\right|\leq C\omega (d(M)).
\end{equation}

\noindent
Let $\bar{\lambda}$ be a unit vector. We have
%(30)
\begin{equation}
\begin{split}
g'_{2\bar{\lambda}}(M) &={{\left( g_2(\widehat{M})-g(\widehat{M}) \right)}^{\prime}_{\bar{\lambda}\mid}}_{\widehat{M}=M}=  \\
&={{\left(\frac{1}{\left| B^{\star}(\widehat{M})\right|}\int\limits_{B^{\star}(\widehat{M})}\left( g_{1}(\widetilde{M})-g(M)\right)\, dm_{3}(\widetilde{M})\right)}'_{\bar{\lambda}\mid}}_{\widehat{M}=M}=  \\
&={{\left(\frac{1}{\left| B^{\star}(\widehat{M})\right|}\right)}'_{\bar{\lambda}\mid}}_{\widehat{M}=M}\int\limits_{B^{\star}(M)} \left( g_1(\widetilde{M})-g(M)\right)\, dm_{3}(\widetilde{M}) +  \\
&+\frac{1}{\left| B^{\star}(M)\right|} {{\left(\int\limits_{B^{\star}(\widehat{M})}\left( g_{1}(\widetilde{M})-g(M)\right)\, dm_{3}(\widetilde{M})\right)}'_{\bar{\lambda}\mid}}_{\widehat{M}=M}=  \\
& =-\frac{{\left| B^{\star}(M)\right|}'_{\bar{\lambda}}}{{\left| B^{\star}(M)\right|}^{2}}\int\limits_{B^{\star}(M)}\left( g_{1}(\widetilde{M})-g(M)\right)\; dm_{3}(\widetilde{M}) +  \\
& + \frac{1}{\left| B^{\star}(M)\right|}\int\limits_{\partial B^{\star}(M)} \left( \bar{n}(\widetilde{M}),\bar{\lambda}\right)\left( g_1(\widetilde{M})-g(M)\right)\, dm_{2}(\widetilde{M}),
\end{split}
\end{equation}

\noindent
where $\bar{n}(\widetilde{M})$ in the last integral is the unit vector of the outer normal to the sphere $\partial B^{\star}(M)$  and $dm_{2}(\widetilde{M})$ denotes the two-dimensional surface measure on the sphere $\partial B^{\star}(M)$.

%page 10
Applying estimates (23) and (19) to   formula (30), we find that

%(31)
\begin{equation}
\left| g'_{2\bar{\lambda}}(M)\right|\leq C\frac{\omega (d(M))}{d(M)},
\end{equation}

\noindent
hence

%(32)
\begin{equation}
\left|\nabla g_{2}(M)\right|\leq C\frac{\omega (d(M))}{d(M)}.
\end{equation}
Repeating the same reasoning as in (30), we obtain by (28), (31), and (32) the following estimate  for $g_0$:
%(33)
\begin{equation}
\left| g'_{0\bar{\lambda}}(M)\right|\leq C\frac{\omega (d(M))}{d(M)}.
\end{equation}

Let  $\bar{\lambda}$ and $\bar{\mu}$ be two arbitrary unit vectors. Then
\begin{eqnarray}
%\begin{split}
&{g''_{0}}_{\bar{\lambda}\bar{\mu}}(M)={{\left( g_0(\widehat{M})-g(M) \right)}''_{\bar{\lambda}\bar{\mu}\mid}}_{\widehat{M}=M}= \nonumber \\
&={{\left(\frac{1}{\left| B^{\star}(\widehat{M})\right|}\int\limits_{B^{\star}(\widehat{M})}\left( g_{2}(\widetilde{M})-g(M)\right)\, dm_{3}(\widetilde{M})\right)}''_{\bar{\lambda}\bar{\mu}\mid}}_{\widehat{M}=M}  \nonumber \\
&={{\left(\frac{1}{\left| B^{\star}(\widehat{M})\right|}\right)}''_{\bar{\lambda}\bar{\mu}\mid}}_{\widehat{M}=M}\int\limits_{B^{\star}(M)} \left( g_2(\widetilde{M})-g(M)\right)\, dm_{3}(\widetilde{M})  \nonumber \\
&+{{\left(\frac{1}{\left| B^{\star}(\widehat{M})\right|}\right)}'_{\bar{\lambda}\mid}}_{\widehat{M}=M}
{{\left(\int\limits_{B^{\star}(\widehat{M})} \left( g_2(\widetilde{M})-g(M)\right)\, dm_{3}(\widetilde{M})\right)}'_{\bar{\mu}\mid}}_{\widehat{M}=M} \nonumber \\
&+{{\left(\frac{1}{\left| B^{\star}(\widehat{M})\right|}\right)}'_{\bar{\mu}\mid}}_{\widehat{M}=M}
{{\left(\int\limits_{B^{\star}(\widehat{M})} \left( g_2(\widetilde{M})-g(M)\right)\, dm_{3}(\widetilde{M})\right)}'_{\bar{\lambda}\mid}}_{\widehat{M}=M} \nonumber \\
&+\frac{1}{\left| B^{\star}(M)\right|} {{\left(\int\limits_{B^{\star}(\widehat{M})}\left( g_{2}(\widetilde{M})-g(M)\right)\, dm_{3}(\widetilde{M})\right)}''_{\bar{\lambda}\bar{\mu}\mid}}_{\widehat{M}=M}  \\
&=-{{\left(\frac{{\left| B^{\star}(\widehat{M})\right|}'_{\bar{\lambda}}}{{\left| B^{\star}(\widehat{M})\right|}^{2}}\right)}'_{\bar{\mu}\mid}}_{\widehat{M}=M}\int\limits_{B^{\star}(M)}\left( g_{2}(\widetilde{M})-g(M)\right)\; dm_{3}(\widetilde{M}) \nonumber \\
&-\frac{{{\left| B^{\star}(\widehat{M})\right|}'_{\bar{\lambda}\mid}}_{\widehat{M}=M}}{{\left| B^{\star}(M)\right|}^{2}} \nonumber
     \int\limits_{\partial B^{\star}(M)}\left(\left( \bar{n}(\widetilde{M}),\bar{\mu}\right)+\left( r^{\star}(M)\right)'_{\bar{\mu}}\right)\cdot\left( g_2(\widetilde{M})-g(M)\right)\, dm_{2}(\widetilde{M}) \nonumber  \\
&-\frac{{{\left| B^{\star}(\widehat{M})\right|}'_{\bar{\mu}\mid}}_{\widehat{M}=M}}{{\left| B^{\star}(M)\right|}^{2}}
     \int\limits_{\partial B^{\star}(M)}\left(\left( \bar{n}(\widetilde{M}),\bar{\lambda}\right)+\left( r^{\star}(M)\right)'_{\bar{\lambda}}\right)\cdot\left( g_2(\widetilde{M})-g(M)\right)\, dm_{2}(\widetilde{M}) \nonumber \\
&+ \frac{1}{\left| B^{\star}(M)\right|}
     {\left( \int\limits_{\partial B^{\star}(\widehat{M})}\left(\left( \bar{n}(\widetilde{M}),\bar{\lambda}\right)+\left( r^{\star}(\widehat{M})\right)'_{\bar{\lambda}}\right)
		\left( g_2(\widetilde{M})-g(M)\right)\, dm_{2}(\widetilde{M})\right)}'_{{\left.\bar{\mu}\right|}_{\widehat{M}=M}} \nonumber
%\end{split}
\end{eqnarray}

%page12
Now we take into account that

%(35)
\begin{equation}
\begin{split}
&{{\left( \int\limits_{\partial B^{\star}(\widehat{M})}\left(\left( \bar{n}(\widetilde{M}),\bar{\lambda}\right)+\left( r^{\star}(\widehat{M})\right)'_{\bar{\lambda}}\right)\left( g_2(\widetilde{M})-g(M)\right)\, dm_{2}(\widetilde{M})\right)}'_{\bar{\mu}\mid}}_{\widehat{M}=M}= \\
&=\int\limits_{\partial B^{\star}(M)}{{\left( {r^{\star}(\widehat{M})}'_{\bar{\lambda}}\right)}'_{\bar{\mu}\mid}}_{\widehat{M}=M}\cdot\left( g_2(\widetilde{M})-g(M)\right)\, dm_{2}(\widetilde{M})+ \\
& + 2\int\limits_{\partial B^{\star}(M)}\frac{{( r^{\star}(M) )}'_{\bar{\lambda}}{( r^{\star}(M) )}'_{\bar{\mu}}}{r^{\star}(M)}\left( g_2(\widetilde{M})-g(M)\right)\, dm_{2}(\widetilde{M})+  \\
&+\int\limits_{\partial B^{\star}(M)}{(r^{\star}(\widehat{M}))}'_{\bar{\lambda}}\cdot {(g_2(\widetilde{M}))}'_{\bar{\mu}}\, dm_{2}(\widetilde{M})+ \\
&+2\int\limits_{\partial B^{\star}(M)} \frac{\left(\bar{n} (\widetilde{M}),\bar{\lambda}\right){(r^{\star}(M))}'_{\bar{\mu}}}{r^{\star}(M)}\left( g_2(\widetilde{M})-g(M)\right)\, dm_{2}(\widetilde{M})+  \\
&+\int\limits_{\partial B^{\star}(M)}\left(\left( \bar{n}(\widetilde{M}),\bar{\lambda}\right)+{\left( r^{\star}(M)\right)}'_{\bar{\lambda}}\right) {(r^{\star}(M))}'_{\bar{\mu}} {(g_2(\widetilde{M}))}'_{\bar{n}(\widetilde{M})}\, dm_{2}(\widetilde{M}).
\end{split}
\end{equation}
%page13
Combining   estimates (23), (24), (28), and (33) and equalities (34) and (35), we find that

\begin{equation*}
\left| {g_0}_{\bar \lambda\bar \mu }(M)\right|\leq C\frac{\omega (d(M))}{d^2(M)},
\end{equation*}

\noindent
which implies

\begin{equation*}
\left|\nabla^2 g_{0}(M)\right|\leq C\frac{\omega (d(M))}{d^2(M)},
\end{equation*}

\noindent
and finally,

%(36)
\begin{equation}
\left|\Delta g_{0}(M)\right|\leq C\frac{\omega (d(M))}{d^2(M)}.
\end{equation}

Inequalities (29), (33), and (36) conclude the proof of   Theorem~4 with a slight change in notation: we have produced a required function $g_0$.

\subsection{Proof of Theorem~5}
\label{sec:3}

%{\bf 4. Proof of Theorem~3b}

Now we  proceed to the proof of Theorem~3b. Consider the sets $\widetilde{\Omega}_{\!n}$ and $\Omega_{\!n}$ defined in (9) and (10). The boundaries of $\Omega_n$ and $\widetilde{\Omega}_{\!n}$ consist of a finite number of subsets of spheres of radii $2\Lambda_{\!n}$ and $\Lambda_{\!n}$; the total area of these spheres is

%page14
%(37)
\begin{equation}
4\pi \left( (2^{n+1}+1)\cdot\Lambda_{\!n}^{2} + (2^n+1)\cdot 4\Lambda_{\!n}^{2}\right)\leq C\cdot 2^n\cdot ({2^{-n}})^2 = C\cdot 2^{-n}
\end{equation}
We fix a point $M_0\in\mathbb{R}^3\setminus L$ and choose $n$ such that $M_0\notin\widetilde{\Omega}_{\!n}$. Assume that $f_0$ is a pseudoharmonic extension of $f$ and that  $R_0$ is chosen so large that $f_0(M)\equiv 0$ outside  the ball $B_{R_0}({O})$ and $M_0\in B_{R_0}({O})$. We denote by $\Sigma_{\!n}$ the connected component of the set $B_{R_0}({O})\setminus \Omega_{\!n}$ containing the point $M_0$. Now we use   the classical formula

%(38)
\begin{equation}
\begin{split}
&f_0(M_0)=\frac{1}{4\pi} \int\limits_{\partial\sum_{n}} \left(f_0(M)\right)'_{\bar{n} (M)} \, \frac{1}{\rho_{M_0}(M)} \, dS(M)- \\
&-\frac{1}{4\pi}\int\limits_{\partial\sum_{n}} \, f_0(M)\left(\frac{1}{\rho_{M_0}(M)}\right)'_{\bar{n} (M)} \, dS(M)- \\
&-\frac{1}{4\pi}\int\limits_{\partial\sum_{n}} \,\frac{\Delta f_0(M)}{\rho_{M_0}(M)} \, dm_{3}(M),
\end{split}
\end{equation}

\noindent
where $\rho_{\!M_0}(M)\MYdef \|M_0 M\|$, $\vec n(M)$ is the outer unit normal at $M\in\partial \Sigma_{\!n}$  to the domain $\Sigma_{\!n}$, $dS(M)$ is the two-dimensional measure on $\partial \Sigma_{\!n}$,  and $m_3$ is the three-dimensional Lebesgue measure in $\mathbb{R}^{\!3}$.

%page15
We take into account that $f_0(M)\equiv 0$ and ${\left( f_0(M)\right)}'_{\pi (M)}\equiv 0$ for $M\in \partial B_{R_0}({O})$. This implies that the integrals in (38) are calculated over the domain $\partial \Sigma_{\!n}\bigcap\partial \widetilde{\Omega}_{\!n}$ whose two-dimensional measure does not exceed $c\cdot 2^{-n}$. The    construction of $\widetilde{\Omega}_{\!n}$ gives the estimates $c'\cdot 2^{-n}\leq d(M)\leq c''\cdot 2^{-n}$, $M\in \partial \widetilde{\Omega}_{\!n}$, with some constants $c',\; c'' > 0$, and   condition (6) yields a sequence $\{\alpha_n\}_{n=1}^{\infty}$, $\alpha_n\to 0$, such that

%(39)

\begin{equation}
\left|\ \left( f_0(M)\right)'_{\bar{n} (M)}\right|\leq C \alpha_{n}\left( d(M)\right)^{-1}, \, M\in \Sigma_n.
\end{equation}
Using (39) and the above argument, we obtain

%(40)

\begin{equation}
\left|\frac{1}{4\pi} \int\limits_{\partial\sum_{n}} \left(f_0(M)\right)'_{\bar{n} (M)} \, \frac{1}{\rho_{M_0}(M)} \, dS(M) \right| \leq C\alpha_{n} \cdot 2^{n}\cdot 2^{-n}=C\alpha_{n}\\
\end{equation}

\noindent
and

\begin{equation}
\left|-\frac{1}{4\pi} \int\limits_{\partial\sum_{n}} f_0(M) \cdot \, \left(\frac{1}{\rho_{M_0}(M)}\right)'_{\bar{n} (M)} \, dS(M) \right| \leq C\cdot 2^{-n}. \end{equation}
Formula (38) and estimates (40) and (41) imply the relation
\begin{equation}
\begin{split}
f_0(M_0)=-\frac{1}{4\pi} \int\limits_{\sum_{n}}  \frac{\Delta f_0(M)}{\rho_{M_0}(M)} \, dm_3(M) + O(\alpha_{n} + 2^{-n})
\end{split}
\end{equation}
Passing to the limit in (42), we get
\begin{equation}
\begin{split}
f_0(M_0)=-\frac{1}{4\pi} \int\limits_{B_{R_0}(\mathbb{O})}  \frac{\Delta f_0(M)}{\rho_{M_0}(M)} \, dm_3(M).
\end{split}
\end{equation}

%page16
We will check below that the integral in   (43) is continuous on $\mathbb{R}^{3}$. Equality (43) and the continuity of both sides of it on $\mathbb{R}^{3}$ allows us to take in (43) an arbitrary point $M_0$ of $\mathbb{R}^{3}$. In particular, we can take $M_0\in L$. Bearing this in mind, we take $M_1, M_2 \in L$, $M_1\neq M_2$ and obtain

%(44)
\begin{equation}
\begin{split}
&f(M_2)-f(M_1)=\frac{1}{4\pi} \int\limits_{B_{R_0}(\mathbb{O})}  \frac{\Delta f_0(M)}{\rho_{M_1}(M)} \, dm_3(M)-\\
&- \frac{1}{4\pi} \int\limits_{B_{R_0}(\mathbb{O})}  \frac{\Delta f_0(M)}{\rho_{M_2}(M)} \, dm_3(M) =\\
&= \frac{1}{4\pi} \int\limits_{B_{2\|M_1M_2\|}(M_1)}  \frac{\Delta f_0(M)}{\rho_{M_1}(M)} \, dm_3(M) -\\
&-\frac{1}{4\pi} \int\limits_{B_{2\|M_1M_2\|}(M_1)}  \frac{\Delta f_0(M)}{\rho_{M_2}(M)} \, dm_3(M)+\\
&+\frac{1}{4\pi} \int\limits_{B_{R_0}(\mathbb{O})\setminus B_{2\|M_1M_2\|}(M_1)}  \left(\frac{1}{\rho_{M_1}(M)}-\frac{1}{\rho_{M_2}(M)}\right) \Delta f_0(M)\, dm_3(M) \\
& \stackrel{\rm def}= I_1-I_2+I_3.
\end{split}
\end{equation}

Using assumptions (6), (7), and (8) of Theorem 3b, we conclude that

%(45) page17
\begin{equation}
\begin{split}
&|I_2|\leq \frac{1}{4\pi} \int\limits_{B_{3\|M_1M_2\|}(M_2)}  \frac{|\Delta f_0(M)|}{\rho_{M_2}(M)} \, dm_3(M)\leq \\
&C\int\limits_{B_{3\|M_1M_2\|}(M_2)}  \frac{\omega(d(M))}{d^2(M) \rho_{M_2}(M)} \, dm_3(M)=\\
&=C\sum_{n=0}^{\infty}\int\limits_{B_{3\cdot 2^{-n}\|M_1M_2\|}(M_2)\setminus B_{3\cdot 2^{-n-1}\|M_1M_2\|}(M_2)}  \frac{\omega(d(M))}{d^2(M) \rho_{M_2}(M)} \, dm_3(M)\leq \\
&C\sum_{n=0}^{\infty}\frac{2^n}{\|M_1M_2\|}\int\limits_{B_{3\cdot 2^{-n}\|M_1M_2\|}(M_2)\setminus B_{3\cdot 2^{-n-1}\|M_1M_2\|}(M_2)}  \frac{\omega(d(M))}{d^2(M)} \, dm_3(M)\leq \\
&C\sum_{n=0}^{\infty}\frac{2^n}{\|M_1M_2\|}\int\limits_{B_{3\cdot 2^{-n}\|M_1M_2\|}(M_2)}  \frac{\omega(d(M))}{d^2(M)} \, dm_3(M)
\end{split}
\end{equation}

Without loss of generality, we may assume that $\|M_1M_2\|\leq \widetilde{C}\|AB\|$ with a constant $\widetilde{C}$ such that   $B_{3\|M_1M_2\|}(M_2)\subset \widetilde{\Omega}_{0}$ for  $M_2\in L$, where $\widetilde{\Omega}_{0}$ is the set defined in (9).
Let $\sigma_{nk}=B_{3\cdot 2^{-n}\|M_1M_2\|} (M_2)\bigcap \Omega_{k}$, where the sets $\Omega_{k}$ are defined in (10). Then we can rewrite  a summand in (45) in the following way:

%(46)
\begin{equation}
\begin{split}
&\int\limits_{B_{3\cdot 2^{-n}\|M_1M_2\|}(M_2)}  \frac{\omega(d(M))}{d^2(M)} \, dm_3(M)= \\
&=\sum_{k=0}^{\infty}\int\limits_{\sigma_{nk}} \frac{\omega(d(M))}{d^2(M)} \, dm_3(M) = \\
&=\sum_{k=k(n)}^{\infty}\int\limits_{\sigma_{kn}}  \frac{\omega(d(M))}{d^2(M)} \, dm_3(M).
\end{split}
\end{equation}

The index $k(n)$ in  (46) means the smallest $k$ such that $\Omega_k\cap B_{3\cdot 2^{-n}\|M_1M_2\|}(M_2)\neq \varnothing$. Inequalities (11) imply the following  important estimates:

%(47)
\begin{equation}
2^{-k(n)}\asymp 2^{-n}\cdot \|M_1M_2\|,
\end{equation}

%(48)
\begin{equation}
d(M)\asymp 2^{-k}, \, M\in \sigma_{kn}
\end{equation}

%page18

Let $\widetilde{\sigma}_{nk}=B_{3\cdot 2^{-n}\|M_1M_2\|} (M_2)\bigcap \widetilde{\Omega}_{k}$, then $\sigma_{nk}\subset \widetilde{\sigma}_{nk}$ and $m_3 \sigma_{nk}\leq m_3 \widetilde{\sigma}_{nk}$. Since \\ $\Omega_{\nu}\bigcap B_{3\cdot 2^{-n}\|M_1M_2\|} (M_2)=\varnothing$ for $\nu < k(n)$, we see that, for $k\ge k(n)$, the center of each ball constituent  of $\widetilde{\sigma}_{nk}$ of radius $2^{-k}$  lies on a subarc of $L$ of  length $\leq C\cdot 2^{-n}\cdot \|M_1M_2\|$, which implies that the number $N_{n,k}$ of such balls does not exceed $C\cdot 2^{-n}\cdot \|M_1M_2\|\cdot 2^k$. Hence

%(49)
\begin{equation}
m_3 \widetilde{\sigma}_{nk}\leq C N_{nk}\cdot 2^{-3k}\leq C\cdot 2^{-n}\cdot\|M_1M_2\|\cdot 2^k \cdot 2^{-3k}=C\cdot 2^{-n-2k}\|M_1M_2\|.
\end{equation}
Finally, combining   estimates (47), (48), and (49), we obtain

%(50)
\begin{equation}
\begin{split}
&\sum_{k=k(n)}^{\infty}\int\limits_{\sigma_{kn}} \frac{\omega(d(M))}{d^2(M)} \, dm_3(M)\leq  \\
&C\sum_{k=k(n)}^{\infty} 2^{2k}\omega(2^{-k}) m_3(\widetilde{\sigma}_{nk})\leq \\
&C\sum_{k=k(n)}^{\infty} 2^{2k}\omega(2^{-k})\cdot 2^{-n-2k}\|M_2M_2\|=\\
&=C 2^{-n}\|M_1M_2\|\cdot\sum_{k=k(n)}^{\infty} \omega(2^{-k}).
\end{split}
\end{equation}

The first assumption in (1) concerning $\omega(t)$ gives the inequality

%(51)
\begin{equation}
\sum_{k=k(n)}^{\infty} \omega(2^{-k})\leq C\omega(2^{-k(n)})\leq C\omega(2^{-n}\cdot \|M_1M_2\|).
\end{equation}

So, formulas (46), (50), and (51) imply the estimate
%(52)
\begin{equation}
\int\limits_{B_{3\cdot 2^{-n}\|M_1M_2\|} (M_2)} \frac{\omega(d(M))}{d^2(M)} \, dm_3(M)\leq C 2^{-n}\cdot\|M_1M_2\|\cdot\omega(2^{-n}\cdot \|M_1M_2\|).
\end{equation}

Let us substitute (52) into (45). Using (1), we obtain

%page19
%(53)
\begin{equation}
\begin{split}
&\sum_{n=0}^{\infty}\frac{2^n}{\|M_1M_2\|}\int\limits_{B_{3\cdot 2^{-n}\|M_1M_2\|} (M_2)} \frac{\omega(d(M))}{d^2(M)} \, dm_3(M)\leq \\ &C\sum_{n=0}^{\infty}\frac{2^n}{\|M_1M_2\|}\cdot 2^{-n}\cdot\|M_1M_2\|\omega(2^{-n}\cdot\|M_1M_2\|)= \\
&=C\sum_{n=0}^{\infty} \omega(2^{-n}\cdot\|M_1M_2\|)\leq C\omega(\|M_1M_2\|),\\
\end{split}
\end{equation}
which means that $|I_2|\leq C\omega(\|M_1M_2\|)$.

The same arguments show that $|I_1|\leq C\omega(\|M_1M_2\|)$. To estimate the term $I_3$, we  use the second part of   assumption (1) concerning the function $\omega(t)$. We notice that, for all $M\notin B_{2\|M_1M_2\|(M_1)}$,
we have the inequality
%(54)
\begin{equation}
\left |\frac{1}{\rho_{M_1}(M)}-\frac{1}{\rho_{M_2}(M)}\right |\leq C\frac{\|M_1M_2\|}{\rho_{M_1}^2(M)}.
\end{equation}
Now, using (54) and (8), we obtain

%(55)
\begin{equation}
\begin{split}
&|I_3|\leq C\int\limits_{B_{R_0}(\mathbb{O})\setminus B_{2\|M_1M_2\|}(M_1)}  \frac{\|M_1M_2\|}{\rho_{M_1}^2(M)} |\Delta f_0(M)|\, dm_3(M)\leq \\
&C\int\limits_{B_{R_0}(\mathbb{O})\setminus B_{2\|M_1M_2\|}(M_1)}  \frac{\|M_1M_2\|}{\rho_{M_1}^2(M)} \frac{\omega(d(M))}{d^2(M)}\, dm_3(M)\leq \\
&C\sum_{n=1}^{\infty}\int\limits_{(B_{2^{n+1}\|M_1M_2\|}(M_1)\setminus B_{2^{n}\|M_1M_2\|}(M_1))\bigcap B_{R_0}(\mathbb{O})}
\frac{\|M_1M_2\|}{\rho_{M_1}^2(M)} \frac{\omega(d(M))}{d^2(M)}\, dm_3(M)\leq \\
&C\sum_{n=1}^{\infty}\frac{\|M_1M_2\|}{2^{2n}\|M_1M_2\|^{2}}\int\limits_{(B_{2^{n+1}\|M_1M_2\|}(M_1)\setminus B_{2^{n}\|M_1M_2\|}(M_1))\bigcap B_{R_0}(\mathbb{O})}\frac{\omega(d(M))}{d^2(M)}\, dm_3(M)\stackrel{\rm def}=\\
&C\frac{1}{\|M_1M_2\|}\sum_{n=1}^{\infty}\frac{1}{2^{2n}} C_n.
\end{split}
\end{equation}

Now, repeating the same reasoning as we used to get (47)--(52), we obtain the estimate
%(56)
\begin{equation}
C_n\leq C\cdot 2^{n}\cdot \|M_1M_2\|\cdot \omega(2^{n}\|M_1M_2\|).
\end{equation}

Combining (55) and (56), we see that

%(57)
\begin{equation}
\begin{split}
&|I_3|\leq C\frac{1}{\|M_1M_2\|}\sum_{n=1}^{\infty}\frac{1}{2^{2n}}\cdot 2^{n}\cdot \|M_1M_2\|\cdot \omega(2^{n}\|M_1M_2\|)=\\
&=C\sum_{n=1}^{\infty}\frac{\omega(2^{n}\|M_1M_2\|)}{2^n}\leq C\omega(\|M_1M_2\|).
\end{split}
\end{equation}
We made use of the second part of   condition (1) in the last inequality in (57). So, we have proved that $|I_1|,\,|I_2|,\,|I_3|\leq C\omega(\|M_1M_2\|)$, which together with Proposition (44) finishes the proof of  Theorem~5.

%page21
\subsection{Proof of Theorem~1}
\label{sec:4}

%{\bf 5. Proof of the Theorem~1}

We start with some geometrical observations. We divide $L$ by the points $$A=M_{0n},\,M_{1n},\,\ldots M_{2^n n}=B$$ as we did in the definitions (9) and (10) of the domains $\widetilde{\Omega}_{n}$ and $\Omega_{n}$. Let $\Lambda_{n}=2^{-n}\cdot |\Lambda|$,   $C_1\ge 1,\, 0\leq k_0\leq 2^n$, and $B[C_1]=B_{C_1\cdot\Lambda_{n}}(M_{k_0 n})$. Let  $P_0,\, P_1\in \partial B[C_1]\bigcap L$ be such that the subarc $L(P_0,\, P_1)$ of $L$ with the endpoints $P_0$ and $P_1$ is the biggest one if $\partial B[C_1]\bigcap L$ contains more than two points. Then we have  $L(P_0,P_1)\leq C_0\cdot 2C_1\cdot \Lambda_{n}$, and there are at most $[2C_0C_1]+2\leq 2(C_0+1)C_1$ subarcs $L(M_{k,2^n},\, M_{k+1,2^n})$ intersecting $L(P_0,P_1)$. Then it is clear that $$m_3(B[C_1]\bigcap \widetilde{\Omega}_{n})\leq 2(C_0+1)C_1\cdot\frac{4}{3}\pi\cdot (2\Lambda_{n})^3=2(C_0+1)C_1\cdot \frac{32}{3}\pi\Lambda_{n}^{3}.$$

The volume of $B[C_1]$ is equal to $\frac{4}{3}\pi\cdot C_1^{3}\cdot\Lambda_{n}^3$. Therefore,  we can choose  $C_1$ such that $m_3(B[C_1])\ge 2m_3(B[C_1]\bigcap\widetilde{\Omega}_{n})$. We introduce the sets $\beta_{kn},\, 0\leq k\leq 2^n$, as follows: $\beta_{0n}=B_{2\Lambda_{n-2}(M_{0n})},\, \beta_{kn}=B_{2\Lambda_{n-2}}(M_{kn})\setminus\bigcup_{\nu=1}^{k-1}B_{2\Lambda_{n-2}}(M_{\nu n})$. We take a constant $C_1$ in such a way that $m_3(B_{C_1\Lambda_{n}}(M_{kn})\setminus\widetilde{\Omega}_{n-2})\ge \frac{1}{2}m_3(B_{C_1\Lambda_{n}}({O}))$.

The above arguments show that we can choose   $C_1$ depending only on $C_0$. Due to estimates (11) we obtain that the inequality $d(M)\ge 2^{-n+1}$ is  valid for all $M\in B_{C_1(\Lambda_{0})}(M_{kn})\setminus\widetilde{\Omega}_{n-2}$. On the other hand, $d(M)\leq C_1 2^{-n}|\Lambda|$.

Now we proceed   to the definition of   $\upsilon_{2^{-n}}(M)$. Using (46)--(52),
 we obtain

%(58)
\begin{equation}
\int\limits_{B_{2\Lambda_{n-2}}(M_{kn})}\frac{\omega(d(M))}{d^2(M)}\, dm_3(M)\leq C\Lambda_{n-2}\omega(\Lambda_{n-2})
\end{equation}

%page22
Inequality (58) and the definition of the set $\beta_{kn}\subset B_{2\Lambda{n-2}}(M_{kn})$ imply

%(59)
\begin{equation}
\int\limits_{\beta_{kn}}\frac{\omega(d(M))}{d^2(M)}\, dm_3(M)\leq C\Lambda_{n-2}\omega(\Lambda_{n-2}).
\end{equation}

Now we  apply Theorem~4 and construct a pseudoharmonic extension $f_0(M)$ of   $f$. Then (8) and (59) give the relation

%(60)
\begin{equation}
\int\limits_{\beta_{kn}}\, \Delta f_0(M)\, dm_3 = C_{kn}\Lambda_{n-2}\omega(\Lambda_{n-2}),
\end{equation}

\noindent
where $|C_{kn}|\leq C$ for all $n$ and $k$, $0\leq k\leq 2^{n-2}$. We denote by $\chi_{kn}$ the characteristic function of the set $B_{C_1\Lambda_{n}}(M_{kn})\setminus \widetilde{\Omega}_{n-2}$ and put

%(61)
\begin{equation}
\phi_{kn}(M)=\gamma_{kn}\, \Lambda^{-2}_{n}\chi_{kn}(M)\omega (\Lambda_{n}),
\end{equation}

\noindent
where $\gamma_{kn}$ satisfies the condition

%(62)
\begin{equation}
\int\limits_{\beta_{kn}}\, \Delta f_0(M)\, dm_3 + \int\limits_{\mathbb{R}^3}\,\phi_{kn}(M)\, dm_3(M)=0.
\end{equation}
Taking into account (60) and (61) and the definition of the constant $C_1$, we obtain that $|\gamma_{kn}|\leq C$, where $C$ is independent of $k$ and $n$. Further, we define

%(63)
\begin{equation}
\Phi_{n}=\sum_{k=0}^{2^{n-2}}\, \phi_{kn}(M).
\end{equation}
%page23
Preserving the notation $\rho_{M_0}(M)=\|M_0 M\|$, we define  the function $\upsilon_{2^{-n}}(M_0)$ as follows:
%(64)
\begin{equation}
\begin{split}
&\upsilon_{2^{-n}}(M_0)=-\frac{1}{4\pi}\int\limits_{\mathbb{R}^3\setminus\widetilde{\Omega}_{n-2}}\,\frac{\Delta f_0(M)}{\rho_{M_0}(M)}\, dm_3(M)+ \\
&+\frac{1}{4\pi}\int\limits_{\mathbb{R}^3}\,\frac{\Phi_{n}(M)}{\rho_{M_0}(M)}\, dm_3(M)
\end{split}
\end{equation}

%\indent
\subsection{Properties of a function $\upsilon_{2^{-n}}(M_{0})$}
\label{sec:5}

%{\bf Properties of a function $\upsilon_{2^{-n}}(M_0)$}.

Inequality (11) applied to the set $\widetilde{\Omega}_{n-2}$ shows that $(\SP \, \Delta f_0)\cap\Omega_{2^{-n+1}}(L)=\varnothing$ and $(\SP \,  \Phi_n)\cap\Omega_{2^{-n+1}}(L)=\varnothing$. By (64) the function $\upsilon_{2^{-n}}$ is harmonic in $\Omega_{2^{-n+1}}(L)$. Assume that $M_0\in L$. Then, using (43) and (64), we get
%(65)

\begin{equation}
\begin{split}
&\upsilon_{2^{-n}}(M_0)-f(M_0)=-\frac{1}{4\pi}\int\limits_{\mathbb{R}^3\setminus\widetilde{\Omega}_{n-2}}\,\frac{\Delta f_0(M)}{\rho_{M_0}(M)}\, dm_3(M)+ \\
&+\frac{1}{4\pi}\int\limits_{\mathbb{R}^3}\,\frac{\Phi_{n}(M)}{\rho_{M_0}(M)}\, dm_3(M)+\frac{1}{4\pi}\int\limits_{\mathbb{R}^3}\,\frac{\Delta f_0(M)}{\rho_{M_0}(M)}\, dm_3(M)=\\
&=\frac{1}{4\pi}\int\limits_{\widetilde{\Omega}_{n-2}}\,\frac{\Delta f_0(M)}{\rho_{M_0}(M)}\, dm_3(M)+\frac{1}{4\pi}\int\limits_{\mathbb{R}^3}\,\frac{\Phi_{n}(M)}{\rho_{M_0}(M)}\, dm_3(M)=\\
&=\sum_{k=0}^{2^{n-2}}\,\left(\frac{1}{4\pi}\int\limits_{\beta_{kn}}\,\frac{\Delta f_0(M)}{\rho_{M_0}(M)}\, dm_3(M)+\frac{1}{4\pi}\int\limits_{\mathbb{R}^3}\,\frac{\phi_{kn}(M)}{\rho_{M_0}(M)}\, dm_3(M)\right).
\end{split}
\end{equation}

Let $M_0$ belong  to the closed subarc $L\left (M_{k_0,n-2},\, M_{k_{0}+1,n-2} \right )$ of $L$ with the endpoints $M_{k_0,n-2}$ and $M_{k_{0}+1,n-2}$. By (62), we get
%page24
%(66)
\begin{equation}
\sum_{k=0}^{2^{n-2}}=\sum_{k=0}^{k_{0}-2}+\sum_{k=k_{0}-1}^{k_{0}+2}+\sum_{k=k_{0}+3}^{2^{n-2}}\stackrel{\rm def}=\Sigma_1+\Sigma_2+\Sigma_3.
\end{equation}

Now in the same way as in (46)--(53), we get the estimates
%(67)
\begin{equation}
\begin{split}
&\left|\frac{1}{4\pi}\int\limits_{\beta_{kn}}\,\frac{\Delta f_0(M)}{\rho_{M_0}(M)}\, dm_3(M)\right| \\
&\leq C \int\limits_{B_{2\Lambda_{n-2}(M_{k,n-2})}}\,\frac{\omega(d(M))}{d^{2}(M)\cdot 2^{-n+1}}\, dm_3(M) \\
&\leq C\omega(2^{-n+2})\leq C\omega(2^{-n})
\end{split}
\end{equation}

\noindent
for $k_0-1\leq k\leq k_0+2$ because $\rho_{M_0}(M)\ge 2^{-n+1}$ for $M_0\in L$ and $M\in\beta_{kn}$.
Moreover, for all $k$, $0\leq k\leq 2^{n-2}$, we have the inequalities
%(68)
\begin{equation}
\begin{split}
&\left|\frac{1}{4\pi}\int\limits_{\mathbb{R}^3}\,\frac{\phi_{kn}(M)}{\rho_{M_0}(M)}\, dm_3(M)\right|=\frac{|\gamma_n|}{4\pi\, \Lambda^{2}_{n}}\omega(\Lambda_n)\int\limits_{B_{C_1\Lambda_n}(M_{kn})\setminus\widetilde{\Omega}_{n-2}}\, \frac{dm_3(M)}{\rho_{M_0}(M)}\leq \\
&\leq \frac{|\gamma_{kn}|}{4\pi\, \Lambda^{2}_{n}} \omega(\Lambda_n)\cdot \frac{1}{2^{-n+1}}\cdot m_3\left( B_{C_1\Lambda_n}(M_{kn})\setminus\widetilde{\Omega}_{n-2}\right)\leq\\
&\leq C\frac{\omega(\Lambda_n)}{\Lambda^{3}_{n}}\cdot \Lambda^{3}_{n}\leq C\omega(2^{-n}).
\end{split}
\end{equation}
\noindent
Relations  (67) and (68) imply that
%(69)
\begin{equation}
|\Sigma_2|\leq C\omega(2^{-n}).
\end{equation}
Let us suppose now that $k\leq k_0-2$ or $k\ge k_0+3$. Then we transform the summands in $\Sigma_1$ or $\Sigma_3$ as follows:

%(70)
%page25

\begin{equation}
\begin{split}
&\frac{1}{4\pi}\int\limits_{\beta_{kn}}\,\frac{\Delta f_0(M)}{\rho_{M_0}(M)}\, dm_3(M)+\frac{1}{4\pi}\int\limits_{\mathbb{R}^3}\frac{\phi_{kn}(M)}{\rho_{M_0}(M)}\, dm_3(M) = \\
&=\frac{1}{4\pi}\int\limits_{\beta_{kn}}\frac{\Delta f_0(M)}{\rho_{M_0}(M_{k,n-2})}\, dm_3(M) + \\
&+\frac{1}{4\pi} \int\limits_{\beta_{kn}}\, \Delta f_0(M)\left( \frac{1}{\rho_{M_0}(M)}-\frac{1}{\rho_{M_0}(M_{k,n-2})}\right)\, dm_3(M)+\\
&+\frac{1}{4\pi}\int\limits_{B_{C_1\Lambda_n}(M_{kn})\setminus\widetilde{\Omega}_{n-2}}\frac{\gamma_n \omega(\Lambda_n)}{\Lambda^{2}_n}\cdot \frac{1}{\rho_{M_0}(M_{k,n-2})}\, dm_3(M)+\\
&+\frac{1}{4\pi}\int\limits_{B_{C_1\Lambda_n}(M_{kn})\setminus\widetilde{\Omega}_{n-2}}\frac{\gamma_n \omega(\Lambda_n)}{\Lambda^{2}_n}\left( \frac{1}{\rho_{M_0}(M)}-\frac{1}{\rho_{M_0}(M_{k,n-2})}\right)\, dm_3(M)=\\
&=\frac{1}{4\pi}\cdot\frac{1}{\rho_{M_0}(M_{k,n-2})}\left(\, \int\limits_{\beta_{kn}}\, \Delta f_0(M)\, dm_3(M) + \int\limits_{\mathbb{R}^3}\, \phi_{kn}(M)\, dm_3(M) \right)+\\
&+\frac{1}{4\pi} \int\limits_{\beta_{kn}}\, \Delta f_0(M)\left( \frac{1}{\rho_{M_0}(M)}-\frac{1}{\rho_{M_0}(M_{k,n-2})}\right)\, dm_3(M)+\\
&+\frac{1}{4\pi}\int\limits_{B_{C_1\Lambda_n}(M_{kn})\setminus\widetilde{\Omega}_{n-2}}\frac{\gamma_n \omega(\Lambda_n)}{\Lambda^{2}_{n}}\left( \frac{1}{\rho_{M_0}(M)}-\frac{1}{\rho_{M_0}(M_{k,n-2})}\right)\, dm_3(M)=\\
&=\frac{1}{4\pi} \int\limits_{\beta_{kn}}\, \Delta f_0(M)\left( \frac{1}{\rho_{M_0}(M)}-\frac{1}{\rho_{M_0}(M_{k,n-2})}\right)\, dm_3(M)+\\
&+\frac{1}{4\pi}\int\limits_{B_{C_1\Lambda_n}(M_{kn})\setminus\widetilde{\Omega}_{n-2}}\, \frac{\gamma_n \omega(\Lambda_n)}{\Lambda^{2}_{n}}\left( \frac{1}{\rho_{M_0}(M)}-\frac{1}{\rho_{M_0}(M_{k,n-2})}\right)\, dm_3(M) = \\
&=A_k+D_k.
\end{split}
\end{equation}
%page26
We take into account that, for the indices $k$ in question and $M\in \beta_{kn}$, we have
%(71)
\begin{equation}
\begin{split}
&\left|\frac{1}{\rho_{M_0}(M)}-\frac{1}{\rho_{M_0}(M_{k,n-2})}\right| =\left|\frac{1}{\|M_0 M\|}-\frac{1}{\|M_0 M_{k,n-2}\|}\right|\leq \\
&\leq\frac{C\Lambda_{n-2}}{{\|M_0 M_{k,n-2}\|}^2} \leq \frac{C\Lambda_{n-2}}{{|k-k_0|}^2 \Lambda^{2}_{n-2}}\leq\frac{C}{\Lambda_n{|k-k_0|}^2}.
\end{split}
\end{equation}
Since $d(M)\ge 2^{-n+1}$ for $M\in B_{C_1\Lambda_{n}}(M_{k,n-2})\setminus\widetilde{\Omega}_{n-2}$,   inequality (71) is also valid for such $M$ with a different $C$ depending on $C_1$ and $C_0$. Thus, due to (71) we get the following bounds for $A_k$ and $D_k$:

%(72)
\begin{equation}
\begin{split}
&\left| A_k\right|\leq C\int\limits_{\beta_{kn}}\,\frac{\omega(d(M))}{d^{2}(M)}\cdot\frac{1}{\Lambda_n {|k-k_0|}^2}\, dm_3(M)\leq\\
&\leq C\omega(2^{-n})\cdot 2^{-n}\cdot \frac{1}{\Lambda_n {|k-k_0|}^2}\leq C\frac{\omega (2^{-n})}{(k-k_0)^2},
\end{split}
\end{equation}

%(73)
\begin{equation}
\begin{split}
&\left| D_k\right|\leq C\int\limits_{B_{C_1\Lambda_n} (M_{k,n-2})\setminus\widetilde{\Omega}_{n-2}}\,\frac{\omega(\Lambda_{n})}{\Lambda^{2}_{n}}\cdot\frac{1}{\Lambda_n (k-k_0)^2}\, dm_3(M)\leq\\
&\leq C\frac{\omega (2^{-n})}{(k-k_0)^2}.
\end{split}
\end{equation}

\noindent
Consequently, (70), (72), and (73) imply
%(74)
\begin{equation}
\begin{split}
\ &|\Sigma_1|+|\Sigma_3|\leq\sum_{k\leq k_0-2 \atop \text{ or } k\ge k_0+3}\, |A_k|+\sum_{k\leq k_0-2 \atop \text{ or } k\ge k_0+2} \, |D_k|\leq\\
\ &\leq C\omega(2^{-n})\sum_{\nu =1}^{\infty}\frac{1}{\nu^{2}}\leq C\omega(2^{-n}).
\end{split}
\end{equation}
Using (65)--(69) and (74), we have
%(75)
\begin{equation}
\left| \upsilon_{2^{-n}}(M_0)-f(M_0)\right|\leq C\omega(2^{-n}).
\end{equation}
%page27
To get the required estimate (2) for any $\delta > 0$, we choose $n$ such that $2^{-n-1}<\delta \leq 2^{-n}$ and put $\upsilon_{\delta}=\upsilon_{2^{-n}}$;   relation (75) is equivalent to (2).

To verify estimate (3), we begin with the case $\delta =2^{-n}$. Let $\upsilon_{2^{-n}}$ be as before and let $M_0\in \Omega_{2^{-n}}(L)$. We have
\begin{equation*}
\begin{split}
&\left( \upsilon_{2^{-n}}(M_0)\right)'_{\bar{\nu}}=\frac{1}{4\pi}\int\limits_{\mathbb{R}^3\setminus\widetilde{\Omega}_{n-2}}\, \frac{\left(\rho_{M_1}(M)\right)'_{\bar{\nu}\mid_{M_1=M_0}}}{\rho^2_{M_0}(M)}\, \Delta f_0(M)\, dm_3(M) -\\
&-\frac{1}{4\pi}\int\limits_{\mathbb{R}^3\setminus\widetilde{\Omega}_{n-2}}\, \frac{\left(\rho_{M_1}(M)\right)'_{\bar{\nu}\mid_{M_1=M_0}}}{\rho^2_{M_0}(M)}\, \Phi_n(M)\, dm_3(M),
\end{split}
\end{equation*}
where
\noindent
$\bar{\nu}$ is an arbitrary unit vector. Then we get
%(76)
\begin{equation}
\begin{split}
&\left|\left( \upsilon_{2^{-n}}(M_0)\right)'_{\bar{\nu}}\right|\leq C\int\limits_{\mathbb{R}^3\setminus\widetilde{\Omega}_{n-2}}\, \frac{\omega(d(M))}{\rho^2_{M_0}(M) d^2(M)}\, dm_3(M)+ \\
&+ C\int\limits_{\mathbb{R}^3\setminus\widetilde{\Omega}_{n-2}}\, \frac{|\Phi_{n}(M)|}{\rho^2_{M_0}(M)}\, dm_3(M) =\\
&= C\int\limits_{\left(\mathbb{R}^3\setminus\widetilde{\Omega}_{n-2}\right)\cap B_{2^{-n+3}|\Lambda|}(M_0)}\, \frac{\omega(d(M))}{\rho^2_{M_0}(M) d^2(M)}\, dm_3(M) + \\
&+ C\sum_{k=1}^{\infty}\int\limits_{\left( B_{2^{-n+k+3}|\Lambda|}(M_0)\setminus B_{2^{-n+k+2}|\Lambda|}(M_0)\right)\setminus\widetilde{\Omega}_{n-2}}\, \frac{\omega(d(M))}{\rho^2_{M_0}(M) d^2(M)}\, dm_3(M)+\\
&+ C\int\limits_{B_{2^{-n}\cdot |\Lambda^3|}(M_0)\setminus\widetilde{\Omega}_{n-2}}\, \frac{|\Phi_{n}(M)|}{\rho^2_{M_0}(M)}\, dm_3(M)+\\
&+ C\sum_{k=1}^{\infty}\int\limits_{\left( B_{2^{-n+k+3}|\Lambda|}(M_0)\setminus B_{2^{-n+k+2}|\Lambda|}(M_0)\right)\setminus\widetilde{\Omega}_{n-2}}\, \frac{|\Phi_{n}(M)|}{\rho^2_{M_0}(M)}\, dm_3(M).
\end{split}
\end{equation}
%page28
Due to (11)  we have  $d(M)\ge 2^{-n+1}$ for $M\notin \widetilde{\Omega}_{n-2}$. Hence  $$\rho_{M_0}(M)=\|M_0 M\|\ge 2^{-n+1}-2^{-n}=2^{-n}$$ for $M_0\in \Omega_{2^{-n}}(L)$. As in (46)--(52), we obtain
%(77)
\begin{equation}
\begin{split}
&\int\limits_{B_{2^{-n+3}\cdot |\Lambda|}(M_0)\setminus\widetilde{\Omega}_{n-2}}\, \frac{\omega (d(M))}{\rho^{2}_{M_0}(M) d^2(M)}\, dm_3(M)\leq\\
&\leq \frac{C}{2^{-2n}}\int\limits_{B_{2^{-n+3}\cdot |\Lambda |}(M_0)}\, \frac{\omega (d(M))}{d^2(M)}\, dm_3(M)\leq\\
&\leq C\cdot 2^{2n}\cdot 2^{-n+3}\cdot |\Lambda |\cdot \omega (2^{-n+3}|\Lambda |)\leq C\cdot 2^{n}\omega (2^{-n}).
\end{split}
\end{equation}
Using the definition (61) of   $\phi_{kn}$ and the definition (63) of   $\Phi_n$, we obtain analogously that
%(78)
\begin{equation}
\begin{split}
&\int\limits_{B_{2^{-n+3}\cdot |\Lambda |}(M_0)\setminus\widetilde{\Omega}_{n-2}}\, \frac{|\Phi_n(M)|}{\rho^2_{M_0}(M)}\, dm_3(M)\leq\\
&\leq C\cdot 2^{2n}\int\limits_{B_{2^{-n+3}\cdot |\Lambda |}(M_0)}\, |\Phi_n(M)|\, dm_3(M)\leq\\
&\leq C\cdot 2^{n}\omega (2^{-n}).
\end{split}
\end{equation}

Using analogs of   (46)--(52) once again, we get the estimates

%(79)
\begin{equation}
\begin{split}
&\sum_{k=1}^{\infty}\int\limits_{\left( B_{2^{-n+k+3}|\Lambda |}(M_0)\setminus B_{2^{-n+k+2}|\Lambda |}(M_0)\right)\setminus\widetilde{\Omega}_{n-2}}\, \frac{\omega (d(M))}{\rho^2_{M_0}(M) d^2(M)}\, dm_3(M)\leq\\
&\leq C\sum_{k=1}^{\infty} 2^{2n-2k}\int\limits_{B_{2^{-n+k+3}\cdot |\Lambda |}(M_0)}\, \frac{\omega (d(M))}{d^2(M)}\, dm_3(M)\leq\\
&\leq C\sum_{k=1}^{\infty} 2^{2n-2k}\cdot 2^{-n+k}\cdot \omega (2^{-n+k})= C\cdot 2^{n}\sum_{k=1}^{\infty} 2^{-k}\omega (2^{-n+k})\leq\\
&\leq C\cdot 2^{n}\cdot \omega (2^{-n}).
\end{split}
\end{equation}
The last inequality in (79) is a consequence of the second part of assumption (1) concerning $\omega (t)$.
%page29
The definition (61) of   $\phi_{kn}$ allows us to deal with the function $|\Phi_n(M)|$ in the same way as with the expression $\frac{\omega (d(M))}{d^2(M)}$, so we get the relation
%(80)
\begin{equation}
\sum_{k=1}^{\infty}\int\limits_{\left( B_{2^{-n+k+3}|\Lambda |}(M_0)\setminus B_{2^{-n+k+2}|\Lambda |}(M_0)\right)\setminus\widetilde{\Omega}_{n-2}}\, \frac{|\Phi_n(M)|}{\rho^2_{M_0}(M)}\, dm_3(M)\leq C\cdot 2^{n}\cdot \omega (2^{-n}).
\end{equation}
similar to (79).
Combining estimates (76)--(80), we come to the inequality
%(81)
\begin{equation}
\left| \left( \upsilon_{2^{-n}}(M_0)\right)'_{\bar{\nu}}\right| \leq C\cdot 2^{n}\cdot \omega (2^{-n}).
\end{equation}
\noindent
This proves statement (3) for $\delta =2^{-n}$ since the constant $C$ in (81) is independent of $\bar{\nu}$. The case of  arbitrary $\delta >0$ is obtained in the same way as in the proof of statement~(2).

\subsection{Proof of Theorem~2}
\label{sec:6}
%{\bf Proof of the Theorem~2}

We put
%(82)
\begin{equation}
f^{\star}_0(x)=\int\limits_{0}^{x} \frac{\omega (t)}{t}\, dt,\quad x\in [0,\,1],
\end{equation}
and $f^{\star}_0(-x)=f^{\star}_0(x)$.
\noindent
Then   condition (1) implies that $f^{\star}_0(x)\leq C'\omega (x)$ and
%(83)
\begin{equation}
f^{\star}_0(x)\ge \int\limits_{\frac{x}{2}}^{x} \frac{\omega (t)}{t}\, dt\ge \omega \left( \frac{x}{2}\right) \log{2} \ge\widetilde{C}'\omega (x), \quad x\in (0,\,1],
\end{equation}
where
\noindent
$\widetilde{C}'>0$ is independent of $x\in (0,\,1]$.
We have
%(84)
\begin{equation}
f^{\star'}_0(x)=\frac{\omega (x)}{x}, \quad x\in (0,\,1].
\end{equation}

\noindent
Relations (1), (82)--(84) imply that $f^{\star}_0\in H^{\omega} ([0,\,1])$ and $f^{\star}_0(x)\asymp \omega (x)$.
%page(30)
We define $f_0(M)\stackrel{\rm def}= f^{\star}_0(x)$ for $M=(x,\,0,\,0)$. For $A>1$ and $0<x<\frac{1}{A}$, we have
\begin{equation*}
f^{\star}_0(A\, x)>\int\limits_{x}^{A\, x}\, \frac{\omega (t)}{t}\, dt \ge \omega (x) \log{A},
\end{equation*}
\noindent
and so

%(85)
\begin{equation}
\omega (x)\leq \frac{1}{\log{A}} f^{\star}_0(A\, x) \leq \frac{\tilde{C}\omega (A\, x)}{\log{A}}
\end{equation}
Suppose there exist a sequence $\{k_\ell\}_{\ell =1}^{\infty}$ for which conditions (4) and (5) are fulfilled with some constants $C'_1$ and $C'_2$. We may assume that $\lambda_{k_\ell}>4$ for all $\ell$. Every function $V'_{k_{\ell} x}(M)$ is harmonic in the domain $\Omega_{\lambda_{k_\ell}\delta_{k_\ell}}([A_0,\, B_0])$, and (5) gives the following estimate:
%(86)
\begin{equation}
\left| V'_{k_{\ell} x}(M)\right|\leq C'_2 \frac{\omega(\delta_{k_\ell})}{\delta_{k_\ell}}, \; M\in \Omega_{\lambda_{k_\ell}\delta_{k_\ell}}([A_0,\, B_0]).
\end{equation}
Let $r_{\ell} = \frac{1}{2}\lambda_{k_\ell}\delta_{k_\ell}$, and $A_{\ell}=\sqrt{\frac{1}{2}\lambda_{k_\ell}}$. We can use the Poisson integral   representation of the function $V'_{k_{\ell} x}$  harmonic in the ball $B_{2r_{\ell}}({O})$,
%page31
%(87)
\begin{equation}
V'_{k_{\ell} x}(M)=\frac{1}{4\pi r_\ell}\int\limits_{\partial B_{r_\ell}({O})}\, V'_{k_{\ell} x}(P)\frac{r^2_{\ell}-\|{O}M\|^2}{\|M P\|^3}\, dm_{2}(P),
\end{equation}
\noindent
where $M\in B_{r_\ell}({O})$, and $dm_{2}(P)$ denotes the two-dimensional Lebesgue measure on the sphere $\partial B_{r_\ell}({O})$. If $M=(x,\,0,\,0)$, $|x|\leq A_{\ell}\delta_{k_\ell}$, then differentiating the integral (87) with respect to $x$ and taking into account   estimate (86), we obtain the inequality
%(88)
\begin{equation}
\left| V''_{k_{\ell} xx}(M)\right|\leq C'_3\cdot \frac{1}{r_\ell}\max_{P\in \partial B_{r_\ell}({O})}\left| V'_{k_{\ell} x}(P)\right|\leq C'_4\frac{\omega (\delta_{k_\ell})}{r_\ell\delta_{k_\ell}}\leq 2C'_4 \frac{\omega (\delta_{k_\ell})}{\lambda{k_\ell})\delta^{2}_{k_\ell}}.
\end{equation}
Let $x_\ell =A_\ell\delta_{k_\ell}$ and $V^{\star}_k(x)=V_k((x,\,0,\,0))$.  Then (88) implies
%(89)
\begin{equation}
\begin{split}
&\left| V^{\star}_{k_\ell}(x_\ell)+V^{\star}_{k_\ell}(-x_\ell)-2\,V^{\star}_{k_\ell}(0)\right|\leq \max_{|x|\leq x_\ell} \left| V''_{k_{\ell}}(x)\right|\cdot x^2_{\ell}\leq\\
&2C'_4 \frac{\omega (\delta_{k_\ell})}{\lambda{k_\ell})\delta^{2}_{k_\ell}}\cdot A^2_{\ell}\delta^2_{k_\ell}=C'_4\omega(\delta_{k_\ell}).
\end{split}
\end{equation}
From inequality (4) and the definition of $f_0$, it follows that
%(90)
\begin{equation}
\begin{split}
&\left| (f^{\star}_0(x_\ell)-V^{\star}_{k_\ell}(x_\ell)+(f^{\star}_0(-x_\ell)-V^{\star}_{k_\ell}(-x_\ell))-2\, (f^{\star}_0(0)-V^{\star}_{k_\ell}(0)))\right|\leq\\
&\leq 4C'_1\,\omega (\delta_{k_\ell}).
\end{split}
\end{equation}
Estimates (89) and (90) put together imply  that
%(91)
\begin{equation}
\left| f^{\star}_0(x_\ell)+f^{\star}_0(-x_\ell)-2\,f^{\star}_0(0)\right|\leq (C'_4+4C'_1)\omega (\delta_{k_\ell}).
\end{equation}

%page32
\noindent
On the other hand, $f^{\star}$ is an even function, so using relations (82) and (83), we get
%(92)
\begin{equation}
\begin{split}
&f^{\star}_0(x_\ell)+f^{\star}_0(-x_\ell)-2\,f^{\star}_0(0)=2(f^{\star}_0(x_\ell)-f^{\star}_0(0))= \\
&= 2f^{\star}_0(x_\ell)\ge 2\widetilde{C'}\omega (x_\ell)=2\widetilde{C'}\omega (A_\ell\delta_{k_\ell}).
\end{split}
\end{equation}
\noindent
From (91) and (92), we obtain the inequality

%(93)
\begin{equation}
2\widetilde{C'}\omega (A_\ell\delta_{k_\ell})\leq (C'_4+4C'_1)\omega (\delta_{k_\ell}).
\end{equation}

\noindent
Since
 $A_{\ell}\longrightarrow\infty$ as ${\ell\rightarrow \infty}$ and  inequality (93) is fulfilled for all $\ell$,   we have a contradiction with  inequality (85).  Theorem~2 is proved.

%\section{Front matter}
%The author names and affiliations could be formatted in two ways:
%\begin{enumerate}[(1)]
%\item Group the authors per affiliation.
%\item Use footnotes to indicate the affiliations.
%\end{enumerate}
%See the front matter of this document for examples. You are recommended to conform your choice to the journal you are submitting to.
%
%\section{Bibliography styles}

%There are various bibliography styles available. You can select the style of your choice in the preamble of this document. These styles are Elsevier styles based on standard styles like Harvard and Vancouver. Please use Bib\TeX\ to generate your bibliography and include DOIs whenever available.

%Here are two sample references: \cite{Feynman1963118,Dirac1953888}.

\section*{References}

\bibliography{ShirokovAlexeevaJAT1}

\end{document}